\documentclass[10pt,reqno]{amsart}
\usepackage{amsfonts,amsmath,amssymb,amsthm}
\usepackage{latexsym}
\usepackage{eucal}
\usepackage[dvips]{graphics}
\usepackage[english]{babel}
\usepackage[latin1]{inputenc}

\numberwithin{equation}{section}

\begin{document}
\title[Stability for a long-short wave interaction.]
{Comment on ``Orbital Stability of Solitary Wave Solutions for an 
Interaction Equation of Short and Long Dispersive Waves"}

\setlength{\baselineskip}{1.3\baselineskip}
\author[Borys Alvarez-Samaniego]{Borys Alvarez-Samaniego}

\email{balvarez@math.uic.edu}



\thanks{{\textit Date}: 07th May 2005.}
\thanks{Research supported by FAPESP/Brazil under grant No. 2003/09593-2.}

\maketitle

{\scriptsize
 \centerline{Department of Mathematics, Statistics, and Computer Science} 
 \centerline{University of Illinois at Chicago - UIC} 
 \centerline{SEO 322, m/c 249, 851 S Morgan St, Chicago, IL, 60607, USA}}

\begin{abstract}
  J. Angulo and J. F. Montenegro (J. Differential Equations 174 (2001), no. 1, 
  181-199) published a paper about nonlinear stability of 
  solitary waves for an interaction system between a long internal wave 
  and a short surface wave in a two layer fluid considering that the fluid 
  depth of the lower layer is sufficiently large in comparison with the 
  wavelength of the internal wave.  In this note, we show that in a critical 
  step during the proof of  Lemma 2.4 in the above mentioned paper, there is a 
  claim used by the authors which fails to be true.  Lemma 2.4 is crucial 
  for the proof of Lemma 2.7, and for the proof of stability 
  in Theorem 2.1 in the paper before mentioned.  

\vspace{0.5cm}
\noindent
{\textit{MSC:}} 35B35; 35Q51; 35Q35

\vspace{0.5cm}  
\noindent
{\textit{Keywords:}} Nonlinear stability; Solitary waves; Hilbert
transform; Dispersive waves
  
\end{abstract}


Throughout this manuscript we keep the terminology and notations of 
\cite{am:am} as much as possible.  We start (see \cite[Lemma 2.4, p. 188]{am:am})  
by remarking that applying the dominated convergence theorem, we get
\begin{equation*}
 \Big{\|}\Big{[}\frac{1}{c}-\mathcal{K}_{\gamma}^{-1}\Big{]} (\phi_0 f) 
 \Big{\|}^2 = \frac{1}{c^2} \int_{\mathbb{R}} |\xi|^2
 |\widehat{K_{\mu}}(\xi) \widehat{\phi_0 f}(\xi)|^2 d\xi 
 \rightarrow 0,
\end{equation*}
as $\gamma \rightarrow 0^-$, where $\mu=\frac{-c}{\gamma}$.   Now, by 
using the last expression, (2.4) and (2.5) in \cite{am:am}, we obtain   
$\|(\mathcal{M}_{\gamma}-\mathcal{M}_0)f\|\rightarrow 0$ as 
$\gamma \rightarrow 0^{-}$, for all $f \in L^2(\mathbb{R})$.

We note that since  $\widehat{K_{\mu}}(\xi)=\frac{1}{|\xi|+\mu}$, 
it follows that the correct expression for the kernel $K_{\mu}$ is 
given by 
\begin{equation*}
 K_{\mu}(x)  =\left\{
  \begin{array}
  [c]{r}
  \sqrt{\frac{2}{\pi}} \int_0^{+\infty}\frac{\cos(\xi x)}
  {\xi+\mu} d\xi =
  \sqrt{\frac{2}{\pi}} \int_0^{+\infty}\frac{\tau e^{-x\tau}}
  {\mu^2+\tau^2}d\tau,\text{ }x>0,\\
  K_{\mu}(-x),\text{ }x<0,
  \end{array}
\right.
\end{equation*}
where the Fourier transform of $g \in L^1(\mathbb{R})$ is 
defined as $\hat g(\xi)=\frac{1}{\sqrt{2\pi}}\int_{\mathbb{R}} 
g(x)e^{-i\xi x}dx$.  
It is not difficult to see, by using the dominated convergence theorem, 
that 
\begin{equation*}
 K_{\mu}'(x)  =\left\{
  \begin{array}
  [c]{r}
  -\sqrt{\frac{2}{\pi}} \int_0^{+\infty}\frac{\tau^2 e^{-x\tau}}
  {\mu^2+\tau^2}d\tau,\text{ }x>0,\\
  -K_{\mu}'(-x),\text{ }x<0.
  \end{array}
\right.
\end{equation*}
Now, it is easy to verify that $K_{\mu}'\notin L^1(\mathbb{R})$.  Therefore, 
it is not possible to  use inequality (2.6) in \cite{am:am} to conclude that 
$\|\mathcal{M}_{\gamma}-\mathcal{M}_0\|_{B(L^2)}\rightarrow 0$ when  
$\gamma \rightarrow 0^{-}$, as it was claimed in \cite{am:am}.   Hence, claim 
(b) in \cite[Lemma 2.4, p. 187]{am:am}, namely 
$\lim_{\gamma \rightarrow 0^-}\hat \delta(\mathcal{L}_{\gamma}, 
\mathcal{L}_0)=0$, remains unproven.

Finally, we remark that the proof of \cite[Lemma 2.7]{am:am} and 
the proof of stability in \cite[Theorem 2.1]{am:am} make use of 
the spectral structure of the operator $\mathcal{L}_{\gamma}$ 
for $\gamma$ negative and near zero.  These spectral properties 
were obtained in \cite[Lemma 2.4]{am:am} by using the fact  that  
$\lim_{\gamma \rightarrow 0^-}\hat \delta(\mathcal{L}_{\gamma},\mathcal{L}_0)=0$, 
and \cite[Theorem A.2]{am:am} (which follows from    
\cite[Theorem IV.3.16]{k:k}).  In consequence, the part regarding the stability 
of solitary-wave solutions of system (1.1) in \cite[Theorem 2.1]{am:am} 
remains still opened in \cite{am:am}. 

\vspace{0.5cm}
\noindent
{\bf{Acknowledgements.}}  The author wants to express his gratitude 
to Professor Jerry L. Bona (UIC) for his hospitality, and many 
stimulating conversations during his stay at UIC.


\medskip

\end{document}